\documentclass[preprint]{imsart}

\RequirePackage[OT1]{fontenc}
\RequirePackage{amsthm,amsmath}
\RequirePackage[numbers]{natbib}
\RequirePackage[colorlinks,citecolor=blue,urlcolor=blue]{hyperref}

\usepackage{amsmath,amsthm,amssymb,amsfonts}
\usepackage{wasysym,oldgerm,mathrsfs,cancel}
\usepackage{enumerate,stmaryrd}
\usepackage{dsfont,graphicx,yfonts,color,epic,eepic}
\usepackage{xfrac}
\usepackage{makeidx,epstopdf}

\pubyear{2018}
\volume{0}
\issue{0}

\startlocaldefs
\numberwithin{equation}{section}

\usepackage{thmtools}

\theoremstyle{plain}
\newtheorem{thm}{Theorem}[section]
\declaretheorem[style=definition,sibling=thm]{definition}
\declaretheorem[style=definition,sibling=thm]{lemma}
\declaretheorem[style=definition,sibling=thm]{remark}

\newenvironment{Proof}{\noindent\textbf{Proof.\ }}{\hspace*{\fill}$\Box$\medskip}
\newenvironment{lemma0}{\noindent\textbf{Lemma.\ }}{\medskip}

\newcommand{\PR}{\mathds{P}}
\newcommand{\EW}{\mathds{E}}

\newcommand{\dx}{\mathrm{d}}
\newcommand{\supp}{\mbox{supp}}

\newcommand{\C}{\mathcal{C}}
\newcommand{\A}{\mathcal{A}}
\newcommand{\dist}{\mathrm{dist}}

\endlocaldefs

\begin{document}

\begin{frontmatter}
\title{Minimax Euclidean Separation Rates for Testing Convex Hypotheses in $\mathbb{R}^d$\thanksref{t}}
\runtitle{Euclidean Separation Rates}
\thankstext{t}{This research has been partially funded by Deutsche Forschungsgemeinschaft (DFG) through grant CRC 1294 ``Data Assimilation — The seamless integration of data and models'' and by the Deutsche Forschungsgemeinschaft (DFG, German Research Foundation) – 314838170, GRK 2297 MathCoRe.}

\begin{aug}
\author{\fnms{Gilles} \snm{Blanchard}\ead[label=e1]{gilles.blanchard@math.uni-potsdam.de}}

\address{Universit{\"a}t Potsdam, Institut f{\"u}r Mathematik\\
Karl-Liebknecht-Stra{\"ss}e 24-25, 14476 Potsdam, Germany\\
\printead{e1}}

\author{\fnms{Alexandra} \snm{Carpentier}\ead[label=e2]{alexandra.carpentier5@gmail.com}}\\
\and
\author{\fnms{Maurilio} \snm{Gutzeit}\thanksref{t3}\ead[label=e3]{gutzeit@ovgu.de}}

\address{OvGU Magdeburg, Institut f{\"u}r Mathematische Stochastik\\
	Universit{\"a}tsplatz 2, 39106 Magdeburg, Germany\\
\printead{e2,e3}}

\thankstext{t3}{First and corresponding author}

\runauthor{G. Blanchard et al.}

\affiliation{Universit{\"a}t Potsdam and OvGU Magdeburg}

\end{aug}

\begin{abstract}
We consider composite-composite testing problems for the expectation in the Gaussian sequence model where the null hypothesis corresponds to a closed convex subset $\C$ of $\mathbb{R}^d$. We adopt a minimax point of view and our primary objective is to describe the smallest Euclidean distance between the null and alternative hypotheses such that there is a test with small total error probability. In particular, we focus on the dependence of this distance on the dimension $d$ and variance $\frac{1}{n}$ giving rise to the minimax separation rate. In this paper we discuss lower and upper bounds on this rate for different smooth and non-smooth choices for $\C$.
\end{abstract}

\begin{keyword}[class=MSC]
\kwd[Primary ]{60K35}
\kwd{60K35}
\kwd[; secondary ]{60K35}
\end{keyword}
%
\setcounter{tocdepth}{2}
\tableofcontents
\end{frontmatter}

\section{Introduction}

In this paper we consider the problem of testing whether a vector $\mu \in \mathbb R^d$ belongs to a closed convex subset $\C$ of $\mathbb R^d$ with $d\in\mathbb{N}$, based on a noisy observation $X$ obtained from the Gaussian sequence model with variance scaling parameter $n\in\mathbb{N}$, i.e.
\begin{align}\label{eq:model}
X = \mu + \frac{1}{\sqrt{n}} \epsilon,\end{align}
where $\epsilon$ is a standard Gaussian vector. More precisely, in an $l_2$ sense, we aim at finding the order of magnitude of the smallest separation distance $\rho>0$ from $\mathcal{C}$ such that the testing problem
\begin{align}\label{eq:testpb}
H_0 : \mu \in \mathcal C~~~\text{vs.}~~~H_\rho :~ \inf_{c \in \mathcal C}  \Vert \mu - c\Vert_2 \geq \rho,
\end{align}
where $\Vert.\Vert_2$ is the Euclidean distance, can be solved in the following sense: For $\eta\in(0,1)$, we can construct a uniformly $\eta-$consistent test $\varphi$ for (\ref{eq:testpb}), i.e.
\begin{align}\label{eq:testunifdef}
\sup_{\mu \in \mathcal{C}} \mathbb \PR_{\mu}(\varphi=1) + \sup_{\mu\in\mathcal{A}_\rho}\mathbb \PR_{\mu}(\varphi=0)\leq \eta,
\end{align}
where $\mathcal{A}_\rho\subset\mathbb{R}^d$ corresponds to $H_\rho$. We write $\rho^*(\C):=\rho^*(\C,d,n,\eta)$ for the minimax-optimal separation distance $\rho$ of this test, i.e.
$$\rho^*(\C)= \inf\{\rho \geq 0 | \exists~\text{test}~\phi : \sup_{\mu \in \mathcal{C}} \mathbb \PR_{\mu}(\varphi=1) + \sup_{\mu\in\mathcal{A}_\rho}\mathbb \PR_{\mu}(\varphi=0)\leq \eta\},$$
i.e.~the smallest distance $\rho$ that enables the existence of a uniformly consistent test of level $\eta$ for the testing problem~\eqref{eq:testpb}. See Section~\ref{s:setting} for a more precise description of the model and relevant quantities. The theory of minimax testing in general has been profiting very much from the seminal work of Ingster and Suslina, see for example their book \cite{Ing03}.\\

An instance of this problem that was extensively studied is signal detection, i.e. the case where $\mathcal C$ is a singleton, see e.g.~\cite{ingster1998minimax} for an extensive survey of this problem and also \cite{Bar02}. From this literature, we can deduce that the minimax-optimal order of $\rho^*$ in this case is
$$\frac{d^{1/4}}{\sqrt{n}}.$$

In its general form however, this problem is a composite-composite testing problem (i.e. neither $\mathcal{C}$ nor $\mathcal{A}_\rho$ is only a singleton). A versatile way of solving such testing problems was introduced in \cite{gayraud2005adaptive}, where the authors combine signal detection ideas with a covering of the null hypothesis, for deriving minimax optimal testing procedures for composite-composite testing problems, provided that the null hypothesis is not too large (i.e.~that its entropy number is not too large, see Assumption (A3) in~\cite{gayraud2005adaptive}). In this case, the authors prove that the minimax-optimal testing separation rate is the same as the signal detection separation rate, namely $\frac{d^{1/4}}{\sqrt{n}}$. This idea can be generalised also to the case where the null hypothesis is ``too large'' (when Assumption (A3) in~\cite{gayraud2005adaptive} is not satisfied); the approach then implies that an upper bound on the minimax rate of separation is the sum of the signal detection rate and the optimal estimation rate in the null hypothesis $\mathcal C$ -- see~\cite{bull2013adaptive} for an illustration of this for a specific convex shape. Using this technique, one finds that the smaller the entropy of $\mathcal C$, the smaller the separation rate.

This idea has the advantage of generality, but is nevertheless sub-optimal in many simple cases. For instance, if $\mathcal C$ is a half-space, the minimax-optimal separation rate is $\frac{1}{\sqrt{n}}$, which is much smaller than the minimax-opimal signal detection rate, even though a half-space has a much larger entropy (it is even infinite) and larger dimension than a single point. See Section~\ref{General} for an extended discussion on this case. This highlights the fact that for such a testing problem, it is in many cases not the entropy, or size, of the null hypothesis that drives the rate, but rather some other properties of the shape of $\mathcal C$.

In order to overcome the limitations of this approach, some other ideas were proposed. A first line of work can be found in \cite{Bar05}, where the authors consider the general testing problem~\eqref{eq:testpb}, but for separation in $\Vert.\Vert_\infty$-norm instead of $\Vert.\Vert_2$-norm. Since any convex  set can be written as a intersection of half-spaces, they rewrite the problem as a multiple testing problem. This approach is quite fruitful, but the $\Vert.\Vert_\infty$-norm results translate in a non-optimal way to $\Vert.\Vert_2$-norm in terms of the dependence on the dimension $d$, particularly for large $d$. A second main direction that was investigated was to consider testing for some \textit{specific} convex shapes, as e.g.~the cone of positive, monotone, or convex functions, see e.g.~\cite{Jud02}, or also balls for some metrics~\cite{Spo99, carpentier2015testing}. These papers exhibit the minimax-optimal separation distance - or near optimal distance, in some cases of~\cite{Jud02} and~\cite{Spo99} - for the specific convex shapes that are considered, namely cones and smoothness balls. The models considered in these works  are different from our model as they consider functional estimation; also, they do not provide results for more general choices of the null hypothesis. In Sections~\ref{General} and~\ref{Dis}, we derive results for our model and shapes related to those of these papers - namely the positive orthant and the Euclidian ball - in order to relate our work with these earlier results. Finally, a last type of results that are related to our problem is the case where the null hypothesis can be parametrised, see e.g.~\cite{comminges2013minimax} where the authors consider shapes that can be parametrised by a quadratic functional. This approach and their results suggest that the smoothness of the shape of $\C$ has an impact on the testing rate.\\

In this paper, we want to take a more general approach toward the testing problem (\ref{eq:testpb}). In Section \ref{General}, we expose the range of possible separation rates by demonstrating that, without any further assumptions on $\mathcal{C}$, the statement
\begin{equation}\exists w(\eta),w'(\eta)\in(0,\infty):~w(\eta)\frac{1}{\sqrt{n}}\leq\rho^\ast(\mathcal{C})\leq w'(\eta)\frac{\sqrt{d}}{\sqrt{n}}\label{range}\end{equation}
is sharp up to $\ln(d)$-factors. After that, in Sections \ref{Smooth} and \ref{Dis}, we investigate the potential of a geometric smoothness property of the boundary of $\mathcal{C}$. Despite its simplicity, this property takes us quite far: In particular, given any separation rate satisfying~\eqref{range}, it allows for constructing a set $\mathcal{C}$ exhibiting this rate up to $\ln(d)$-factors.

\section{Setting}\label{s:setting}
Let $d,n\in\mathbb{N}$. We consider the $d$-dimensional statistical model
\begin{equation}X=\mu+\frac{1}{\sqrt{n}}\epsilon,\label{Model}\end{equation}
where $\mu\in\mathbb{R}^d$ is unknown and $\epsilon$ is a standard Gaussian vector, written $\epsilon\sim\mathcal{N}(\mathds{O}_d,\mathds{I}_d)$. For $k\in\mathbb{N}$, $\mathds{O}_k$ denotes the origin of $\mathbb{R}^k$ and $\mathds{I}_k\in\mathbb{R}^{k\times k}$ the identity matrix. Clearly, by construction, the variance scaling parameter $n$ may also be interpreted as sample size since the distribution of $X$ is precisely the distribution of the mean of $n$ iid observations from $\mathcal{N}(\mu,\mathds{I}_d)$.

Now, let $\C\subsetneq\mathbb{R}^d$ be closed, nonempty and convex. For $x\in\mathbb{R}^d$ we write
$$\mathrm{dist}(x,\C):=\inf_{c\in\C}\Vert x-c\Vert,$$
where $\Vert \cdot\Vert:=\Vert\cdot\Vert_2$ denotes Euclidean ($l_2$) norm, i.e. $\Vert x\Vert=\sqrt{\sum_{i=1}^k x_i^2}$ for $x\in\mathbb{R}^k$, $k\in\mathbb{N}$. A corresponding open Euclidean ball with center $z\in\mathbb{R}^k$ and radius $r>0$ is denoted $B_k(z,r)$; moreover, we indicate vector concatenation by $[\cdot,\cdot]$, so that, for instance, $[z,1]\in\mathbb{R}^{d}$ if $z\in\mathbb{R}^{d-1}$.\\

Given $\rho>0$, we are interested in the testing problem
\begin{equation}H_0:\ \mu\in\C\ \ \mathrm{vs.}\ \ H_{\rho}:\ \mathrm{dist}(\mu,\C)\geq\rho\label{TP}\end{equation}
and we write $\mathcal{A}_\rho:=\{x\in\mathbb{R}^d\ |\ \mathrm{dist}(x,\C)\geq\rho\}$.
Our goal is to find the smallest value of $\rho$ such that testing $(\ref{TP})$ with prescribed total error probability is possible in a minimax sense, i.e. the quantity
\begin{align*}
\rho^\ast(\C)&:=\rho^\ast(\C,d,n,\eta)\\&~=\inf\{\rho>0\ |\ \exists\ \mathrm{test}\ \varphi:\ \sup_{\mu\in\mathcal{C}}\PR_{\mu}(\varphi=1)+\sup_{\mu\in\mathcal{A}_\rho}\PR_{\mu}(\varphi=0)\leq \eta\}\end{align*}
for some fixed $\eta\in(0,1)$. Here, a test $\varphi$ is a measurable function $\varphi:\mathbb{R}\rightarrow\{0,1\}$.\\

In particular, we focus on the dependence of $\rho^\ast(\C)$ on the dimension $d$ and $n$.  In terms of notation, this is done by using the symbols $\lesssim,\gtrsim$ and $\eqsim$ as follows: For some function $g_{\C}$ that may only depend on $n$ and $d$, we define
$$\rho^\ast(\C)\lesssim g_{\C}(d,n)~~:\Longleftrightarrow~~\exists w(\eta)\in(0,\infty):~\rho^\ast(\C)\leq w(\eta) g_{\C}(d,n).$$
We define in a similar way the symbol $\gtrsim$ (other direction). Finally, if $g_{\C}(d,n)\lesssim \rho^\ast(\C)\lesssim g_{\C}(d,n)$, we write $\rho^\ast(\C)\eqsim g_{\C}(d,n)$; $g_{\C}(d,n)$ then exhibits the \textit{minimax Euclidean separation rate} for (\ref{TP}) or simply \textit{separation rate}.

\begin{remark}
In the proofs for upper bounds on $\rho^\ast(\C)$ it is necessary to consider the type-I and type-II errors $\sup\limits_{\mu\in\mathcal{C}}\PR_{\mu}(\varphi=1)$ and $\sup\limits_{\mu\in\mathcal{A}_\rho}\PR_{\mu}(\varphi=0)$ separately leading to parameters $\alpha,\beta\in(0,\frac{1}{2})$ rather than $\eta$. However, this does not affect the separation rate. For the sake of consistency in notation, we will state the exact constants $w(\eta)$ in upper bounds with $\alpha=\beta=\frac{\eta}{2}$. In these statements and in the proofs, we use the abbreviation $v_x:=\ln\left(\frac{1}{x}\right)$, $x>0$.
\end{remark}


\section{A General Guarantee and Extreme Cases}\label{General}

The quantity $\rho^\ast(\C)$ clearly depends on $\mathcal{C}$.\\
Let us firstly examine a simple, essentially one-dimensional case, namely a half-space.
\begin{thm}\label{HS}
Let $\C=\C_{\mathrm{HS}}:=\mathbb{R}^{d-1}\times (-\infty,0]$ (if $d=1$, $\C_{\mathrm{HS}}=(-\infty,0]$). Then, in the testing problem (\ref{TP}), we have
$$\sqrt{\frac{1}{n}\ln(1+4(1-\eta)^2)}\ \leq\ \rho^\ast(\C)\ \leq\ \sqrt{\frac{8}{n}v_\eta}$$
and therefore
$$\rho^\ast(\C)\ \eqsim\ \frac{1}{\sqrt{n}}.$$
\end{thm}

\begin{remark}\label{GenLow}
As can be seen in the proof (section \ref{HSproof}), this testing problem is essentially equivalent to the problem $\mu=0$ vs. $\mu=\rho$ in dimension $d=1$, so that, alternatively, the rate $\frac{1}{\sqrt{n}}$ can be obtained by analysing the optimal test in the sense of Neyman-Pearson. Furthermore, and in fact closely related to that, note that the lower bound in the previous theorem is valid for any choice of closed convex set $\mathcal{C}$ such that $\mathcal{C}$ and $\mathbb{R}^d\backslash \mathcal{C}$ are non-empty:
$$\frac{1}{\sqrt{n}} \lesssim \rho^\ast(\C).$$
Indeed, we find this rate by considering a fixed pair of points $(\mu_0,\mu_1)\in\C\times\mathcal{A}_{\rho}$ that minimises the distance between $\C$ and $\A_{\rho}$, i.e. $\Vert \mu_0-\mu_1\Vert=\rho$. That seems to have firstly been discussed in \cite{Bur79}; other related (classical) literature would be for instance \cite{Che52} and \cite{Ing00}.\\

\end{remark}
Now, on the other hand, making no additional assumptions about $\C$, a natural choice $\varphi$ for solving (\ref{TP}) is a plug-in test based on confidence balls. This gives rise to the following general upper bound:

\begin{thm}~\label{thm:gen}
Let $\mathcal{C}$ be an arbitrary closed convex subset of $\mathbb{R}^d$ such that $\C$ and $\mathbb R^d \setminus \C$ are non-empty. Then, in the testing problem (\ref{TP}), we have
$$\rho^\ast(\C)\ \leq \ 2\sqrt{\frac{d}{n}+\frac{2}{n}\sqrt{dv_{\eta/2}}+\frac{2}{n}v_{\eta/2}}$$
and therefore
$$\rho^\ast(\C)\ \lesssim\ \frac{\sqrt{d}}{\sqrt{n}}.$$
\label{GenUpp}
\end{thm}

\begin{remark}
Note that this upper bound is the rate of estimation of $\mu$ in $l_2$ norm in the model $(\ref{Model})$ (See Equation (\ref{ConcIn}) in Section \ref{UBGen}). 
\end{remark}

\begin{remark}\label{Genboth}
From Remark~\ref{GenLow} and Theorem~\ref{thm:gen} it is clear that
$$\frac{1}{\sqrt{n}} \lesssim \rho^\ast(\C)\ \lesssim\ \frac{\sqrt{d}}{\sqrt{n}}$$
whenever $\mathcal{C}$ is a closed convex subset of $\mathbb{R}^d$ such that $\C$ and $\mathbb R^d \setminus \C$ are non-empty.\\

\end{remark}
Given this observation, it is natural to ask if the upper bound in Theorem $\ref{GenUpp}$ is also sharp in the sense that there is a choice of $\C$ that requires the separation rate $\frac{\sqrt{d}}{\sqrt{n}}$, at least up to logarithmic factors. It turns out that the answer is yes when $\C$ is taken to be an orthant:

\begin{thm}\label{O}
Let $\C=\C_{\mathrm{O}}:=(-\infty,0]^d$, $d\geq 42$, $\eta\in(0,\frac{8}{9})$ and 
$$M_{\eta}:=\max\left(32,\left\lceil \frac{2}{1-\ln(2)}\ln(d)+1+\frac{2}{1-\ln(2)}\ln\left(\frac{1.8}{\frac{8}{9}-\eta}\right)\right\rceil\right).$$ Then, for the testing problem (\ref{TP}), we have
$$\rho^\ast(\C)\ \geq\ \frac{1}{28}\frac{1}{M_\eta^{3/2}}\frac{\sqrt{d}}{\sqrt{n}}$$
and therefore, if $d$ is large enough in the sense that $M_\eta\leq C\ln(d)$ for some $C>0$,
$$\frac{\sqrt{d}}{\ln(d)^{3/2}\sqrt{n}}\ \lesssim\ \rho^\ast(\C)\ \lesssim\ \frac{\sqrt{d}}{\sqrt{n}}.$$
\label{JNthm}
\end{thm}

This result heavily relies on tailoring the priors such that they have a certain number of moments in common. A related application of this approach can be found in the proof theorem 1 in \cite{Jud02}, see also for instance \cite{Cai11}.

\section{A Simple Smoothness-Type Property}\label{Smooth}
Clearly, the two extreme cases $\C_{\mathrm{HS}}$ and $\C_{\mathrm{O}}$ differ significantly with respect to smoothness of their boundaries. Based on this observation, in order to be able to handle $\rho^\ast(\C)$ more flexibly, we propose to describe convex sets by their boundaries' degree of smoothness, where the boundary of a set $S\in\mathbb{R}^d$ is denoted by $\partial S$ and its closure by $\overline{S}=S\cup\partial S$. To begin with, we examine the potential of the following very simple and purely geometric smoothness concept:
\begin{definition}\label{Rrounded}
Let $R\geq 0$ and $S\subseteq\mathbb{R}^d$ with non-empty interior. $S$ is called \textbf{$\boldsymbol{R}$-rounded} if
\begin{equation}\forall x\in\partial S\ \exists z\in S:\ x\in\overline{B_d(z,R)}\subseteq \overline{S}.\label{RSm}\end{equation}
\end{definition}

\begin{remark}
Note that $R$-rounding is a stronger requirement the higher the value of $R$, i.e. intuitively the degree of the boundary's smoothness grows with $R$. In particular, a half space $\mathcal{C}_{\mathrm{HS}}$ is $\infty$-rounded, a ball $B_d(z,R)$  (with $z\in\mathbb{R}^d$, $R\in(0,\infty)$) is $R$-rounded and the orthant $\mathcal{C}_{\mathrm{O}}$ is $0-$rounded. The definition of $R-$rounding is closely related to the so-called $R$-rolling condition employed in \cite{Ari15}. In fact, $R$-rounding of $S$ is equivalent to saying that $\mathbb{R}^d\backslash S$ fulfils the $R$-rolling condition.

Another related concept worth mentioning is the radius of curvature, though the connection is more subtle: The radius of curvature at a point $x\in\partial S$ would be the radius $r$ of the ball $B$ that best fits $\partial S$ in the sense of a common tangential hyperplane of $\partial S$ and $B$ at $x$ and common analytical curvature, see for instance \cite{Cas96}. Hence, it is possible that the infimum $R$ of these radii $r$ with respect to $x\in\partial S$ corresponds to the parameter $R$ in our previous definition. However, we can then still not easily guarantee that the resulting balls $B$ of the form $B_d(z,R)$ fulfil $\overline{B_d(z,R)}\subseteq\overline{S}$ as required in Definition \ref{Rrounded}.\\

\end{remark}
Since smoothness is usually defined as a local property of a function, we provide a suggestion for how to cast the above concept in that context for a closed convex set $\C$: Given any $x\in\partial\C$, without loss of generality (w.l.o.g.) apply a rotation and translation $G$ such that $x'=G(x)=\mathds{O}_d$ and $\C':=G(\C)\subseteq\mathbb{R}^{d-1}\times[0,\infty)$. Now, assume that there is an $r\in(0,R]$ and a function $f:B\rightarrow[0,\infty)$, where $B=B_{d-1}(\mathds{O}_{d-1},r)$, such that its graph is contained in $\partial\C'$ and $\C'\cap(B\times[0,\infty))$ is contained in the epigraph of $f$ -- see the figure below for an illustration. The following lemma states sufficient conditions for $R$-rounding locally at $x\in\C'$, i.e. at $G^{-1}(x)\in\partial\C$: 

\begin{lemma}\label{Tay}In the situation described in the latter paragraph, if $f$ is twice differentiable on $B$ (i.e. the gradient $\nabla f(\cdot)$ and Hessian matrix $Hf(\cdot)$ exist), the following conditions are sufficient in order that the graph of $f$ remains below $B_d([\mathds{O}_{d-1},R],R)$, i.e. $\mathcal{C}'$ is locally $R$-rounded at $\mathds{O}_d$.  
\renewcommand\arraystretch{2}
$$\left\{\begin{array}{ll}\nabla f(\mathds{O}_{d-1})=\mathds{O}_{d-1},\\
\forall x\in B\backslash\{\mathds{O}_{d-1}\}:\ \ 0\leq \lambda_{\mbox{min}}(Hf(x)),\ \lambda_{\mbox{max}}(Hf(x))\leq\dfrac{1}{R},\end{array}\right.$$
where $\lambda_{\mbox{min}}(\cdot)$ and $\lambda_{\mbox{max}}(\cdot)$ are the lowest and highest eigenvalues of a real symmetric matrix, respectively.\\
\begin{figure*}[!ht]
\centerline{\includegraphics[scale=0.5,trim=0 4cm 0 4cm]{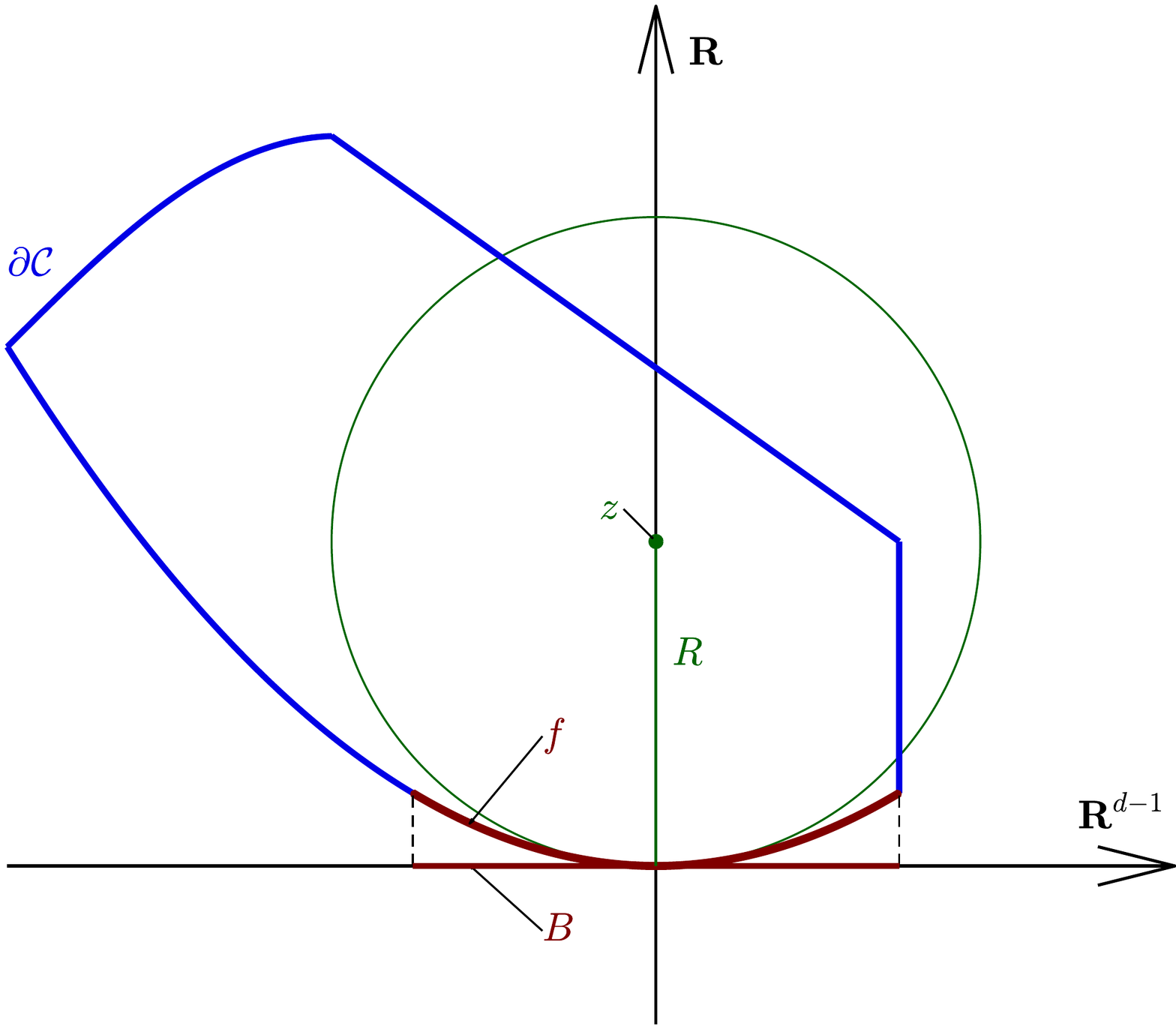}}
\caption{\normalfont In this example $\C$ is only $0$-rounded in the sense of Definition \ref{Rrounded}, but in the local sense of Lemma \ref{Tay}, there are points $x\in\partial\C$ with $0$-rounding, $\infty$-rounding and ``non-degenerate'' rounding such as $\mathds{O}_d$, where, however, the maximum admissible radius $r$ of $B$ is strictly smaller than $R$.}	
\end{figure*}
\end{lemma}\renewcommand\arraystretch{1}
Now, let us examine how the additional assumption of $R-$rounding may affect the general upper bound of Theorem \ref{GenUpp}:
\begin{thm}\label{GenRSm}
If $\C$ is $R-$rounded for some $R\in(0,\infty)$, for the testing problem (\ref{TP}), we have
$$\rho^\ast(\C)\ \leq\ \sqrt{\tfrac{2}{n}v_{\eta/8}}+\tfrac{d}{2nR}+\tfrac{2}{nR}\sqrt{dv_{\eta/4}}+\tfrac{1}{nR}v_{\eta/4}+\sqrt{\tfrac{2}{n}v_{\eta/2}}$$
and therefore, taking Theorem \ref{GenUpp} into account, 
$$\rho^\ast(\C)\ \lesssim\ \max\left(\frac{1}{\sqrt{n}},\min\left(\frac{\sqrt{d}}{\sqrt{n}},\frac{d}{nR}\right)\right)$$
\end{thm}

The following result confirms that this upper bound can be sharp up to $\ln(d)$ factors, namely in the case where $\C$ is taken as an $R$-inflated orthant:

\begin{thm}\label{IO}
Let $d\geq 43$, $\eta\in(0,\frac{8}{9})$ and 
$$\C=\C_{\mathrm{IO}}= \C_{\mathrm{O}} + B(\mathds{O}_d,R)=\bigcup_{z\in\C_{\mathrm{O}}}B(z,R),$$
where $\C_{\mathrm{O}}=(-\infty,0]^d$ is the orthant from Theorem \ref{O}. Furthermore, let
$$M_{\eta}:=\max\left(32,\left\lceil \frac{2}{1-\ln(2)}\ln(d-1)+1+\frac{2}{1-\ln(2)}\ln\left(\frac{1.8}{\frac{8}{9}-\eta}\right)\right\rceil\right).$$
Then, in the testing problem (\ref{TP}), we have with $s=\frac{\sqrt{3}}{28}\frac{1}{M_\eta^{3/2}\sqrt{n}}$
$$\rho^\ast(\C)\ \geq\ \frac{1}{12}\min\left(\frac{(d-1)s^2}{R},\sqrt{3}\sqrt{d-1}\ s\right)$$
and therefore, if $d$ is large enough in the sense that $M_\eta\leq C\ln(d-1)$ for some $C>0$,
$$\rho^\ast(\C)\ \gtrsim\ \max\left(\frac{1}{\sqrt{n}},\min\left(\frac{1}{\ln(d)^3}\cdot \frac{d}{nR},\frac{1}{\ln(d)^{3/2}}\cdot \frac{\sqrt{d}}{\sqrt{n}}\right)\right).$$
\end{thm}

\section{Discussion}\label{Dis}
The concept of $R$-rounding allows for the construction of hypotheses $\C$ with any separation rate $\frac{1}{\sqrt{n}}\lesssim\rho^\ast(\C)\lesssim\frac{\sqrt{d}}{\sqrt{n}}$, up to $\ln(d)$-factors. On the other hand, we must acknowledge that $R-$rounding is too weak a concept to fully describe the difficulty of testing an arbitrary $\C$; an examination of the natural $R$-rounded set, namely a ball of radius $R$, provides clear evidence of this drawback. The result is a direct generalisation of the known rate $\rho^\ast(\C)\sim\frac{d^{\frac{1}{4}}}{\sqrt{n}}$ in the signal detection setting, see \cite{Bar02}.

\begin{thm}\label{Ball}
Let $\eta\in(0,1)$ and $d\geq \ln(2/\eta)$. If $\C=\C_{\mathrm{B}}=B_d(z,R)$, $z\in\mathbb{R}^d$ and $R>0$, for the testing problem (\ref{TP}), we have for $s:=\frac{\sqrt{d-1}}{n}\sqrt{\frac{2}{e}\ln(1+4(1-\eta)^2)}$
$$\rho^\ast(\C)\ \geq\ \frac{s}{2\sqrt{s+R^2}}\ \gtrsim\ \min\left(\frac{d^{\frac{1}{4}}}{\sqrt{n}},\frac{\sqrt{d}}{nR}\right)$$
and also
\begin{eqnarray*}\rho^\ast(\C)&\leq&\min\left(2\sqrt{2}\frac{d^{\frac{1}{4}}}{n^{\frac{1}{2}}}\sqrt{v_{\eta/2}}+3\sqrt{\frac{2}{n}v_{\eta/2}},\frac{2\sqrt{d}}{nR+2\sqrt{nv_{\eta/2}}}\sqrt{v_{\eta/2}}+3\sqrt{\frac{2}{n}v_{\eta/2}}\right)\\
&\lesssim& \max\left(\frac{1}{\sqrt{n}},\min\left(\frac{d^{\frac{1}{4}}}{\sqrt{n}},\frac{\sqrt{d}}{nR}\right)\right).
\end{eqnarray*}
Therefore,
\begin{equation}\rho^\ast(\C)\ \eqsim\ \max\left(\frac{1}{\sqrt{n}},\min\left(\frac{d^{\frac{1}{4}}}{\sqrt{n}},\frac{\sqrt{d}}{nR}\right)\right).\label{BallRate}\end{equation}
\end{thm}

Clearly, Theorem \ref{IO} does not capture this case. As a consequence, future work will be concerned with finding a stronger concept, possibly a localised version of $R-$rounding, that ideally allows for describing $\rho^\ast(\C)$ for any choice of $\C$. However, we suspect this to be quite an ambitious goal. 

\section{Proofs}
\subsection{General Preparations}
\subsubsection{Techniques for Obtaining Lower Bounds}\label{LBGen}
We employ a classical Bayesian approach for proving lower bounds, see references in \cite{Bar02} for its origins. We briefly give the main theoretical ingredients of this approach for our setting:\\

Let $\nu_0$ be a distribution with $S_0:=\supp(\nu_0)\subseteq\C$ and $\nu_\rho$ be a distribution with $S_\rho:=\supp(\nu_\rho)\subseteq\{\mu\in \mathbb{R}^d\ |\ \dist(\mu,\C)=\rho\}$ (priors). For instance, Dirac priors on some $x\in\mathbb{R}$ will be denoted $\delta_x$. Furthermore, for $i\in\{0,\rho\}$, let $\PR_{\nu_i}$ be the resulting distribution of $X$ given $\mu\sim\nu_i$. Now, we see that for any test $\varphi=\mathds{1}_A$, $A\in\mathcal{B}(\mathbb{R}^d)$,
\begin{eqnarray*}
\sup_{\mu\in\mathcal{C}}\PR_{\mu}(\varphi=1)+\sup_{\mu\in\mathcal{A}_\rho}\PR_{\mu}(\varphi=0)&\geq&\PR_{\nu_0}(\varphi=1)+\PR_{\nu_\rho}(\varphi=0)\\
&\geq&1-\frac{1}{2}\Vert \PR_{\nu_\rho}-\PR_{\nu_0}\Vert_{\mathrm{TV}}.\\
&\geq&1-\frac{1}{2}\left(\int_{\mathbb{R}^d}\left(\frac{\dx \PR_{\nu_\rho}}{\dx \PR_{\nu_0}}\right)^2\ \dx \PR_{\nu_0}-1\right)^{\frac{1}{2}},\end{eqnarray*}
see for instance \cite{Bar02}. This justifies the following reasoning used for each lower bound proof in the present paper:\\
Let $\eta\in(0,1)$. For any $\widetilde{\rho}>0$ such that either
$$\frac{1}{2}\Vert \PR_{\nu_0}-\PR_{\nu_{\widetilde{\rho}}}\Vert_{\mathrm{TV}}\leq 1-\eta$$
or
\begin{equation}\int_{\mathbb{R}^d}\left(\frac{\dx \PR_{\nu_\rho}}{\dx \PR_{\nu_0}}\right)^2\ \dx \PR_{\nu_0}\leq 1+4(1-\eta)^2,\label{ChiLB}\end{equation}
it holds that
$$\sup_{\mu\in\mathcal{C}}\PR_{\mu}(\varphi=1)+\sup_{\mu\in\mathcal{A}_\rho}\PR_{\mu}(\varphi=0)\geq\eta$$
and thus, for the testing problem (\ref{TP}), we have
$$\rho^\ast(\C)\geq\widetilde{\rho}.$$

\subsubsection{Concentration Properties of Gaussian and $\chi^2$ Random Variables}\label{UBGen}

We will repeatedly make use of the following classical properties of $N\sim\mathcal{N}(0,\sigma^2)$ and $Z\sim\chi^2_\lambda(d)$ (that is, a $\chi^2$- distribution with $d$ degrees of freedom and noncentrality parameter $\lambda\geq 0$): For any $\delta\in(0,1)$, the following concentration inequalities hold:
\begin{equation}\begin{array}{ccrl}
(\mathrm{I})&&\PR(N\geq \sigma\sqrt{2v_\delta})&\leq\delta\\
(\mathrm{II})&&\PR(Z\geq d+\lambda+2\sqrt{(d+2\lambda)v_\delta}+2v_\delta)&\leq\delta,\\
(\mathrm{III})&&\PR(Z\leq d+\lambda-2\sqrt{(d+2\lambda)v_\delta)}&\leq\delta,\\
\end{array}
\label{ConcIn}
\end{equation}
See \cite{Bir01} for proofs of (\ref{ConcIn}.II) and (\ref{ConcIn}.III). 
\subsubsection{Frequently used Bounds for Expressions Containing Square Roots}
We will employ the following bounds on several occasions which makes it convenient to mention them here.\\
\begin{lemma0}For any $a>0,b\in\mathbb{R}$, we have
\begin{equation}\frac{a}{2\sqrt{a+b^2}}\leq \sqrt{a+b^2}-b\leq\frac{a}{2b}\label{SqrtIn1}\end{equation}
and for any $b> 0$, $a\leq b^2$ we have
\begin{equation}b-\sqrt{b^2-a}\geq\frac{a}{2b}\label{SqrtIn2}\end{equation}
\end{lemma0}\\
\begin{Proof}
	Firstly, through Taylor expansion of $\sqrt{a+b^2}-b$ as a function in $a$, we see that there is a $\xi\in(0,a)$ such that
	$$\sqrt{a+b^2}-b=\frac{a}{2\sqrt{\xi+b^2}}.$$
	Now, with $\xi\geq 0$ and $\xi\leq a$ we obtain the upper and lower bounds in \eqref{SqrtIn1}, respectively. Secondly, explicit calculation tells us that
	$$b-\sqrt{b^2-a}\geq\frac{a}{2b}~\Leftrightarrow~\frac{a^2}{4b^2}\geq 0,$$
	which concludes the proof.
	\end{Proof}

\subsection{Proofs for Section \ref{General}}

\subsubsection{Proof of Theorem \ref{HS}}\label{HSproof}
\begin{Proof}
We prove independently that the order of $\rho^\ast(\C)$ is lower and upper bounded by $\frac{1}{\sqrt{n}}$.
\begin{enumerate}
\item Lower Bound.\\
In accordance with the framework in Section \ref{LBGen}, we verify that the bound holds in the special case $\nu_0=\delta_{\mathds{O}_d}$ and $\nu_{\rho}=\delta_{\rho\cdot e_d}$, where $e_d$ is the last standard basis vector $e_d=[\mathds{O}_{d-1},1]$. Since both the null and alternative hypotheses are simple, the corresponding density functions $F_{\nu_0}(x)$ and $F_{\nu_\rho}(x)$ are readily given and we obtain
\begin{align*}\int_{\mathbb{R}^d}\frac{F_{\nu_\rho}^2}{F_{\nu_0}}(x)\ \dx x&=\sqrt{\frac{n}{2\pi}}^{d}\int_{\mathbb{R}^d}\exp\left(-n(x_d-\rho)^2+\frac{n}{2}x_d^2-\frac{n}{2}(\Vert x\Vert^2-x_d^2)\right)\ \dx x\\
&=\sqrt{\frac{n}{2\pi}}\int_{\mathbb{R}}\exp\left(-n(x_d-\rho)^2+\frac{n}{2}x_d^2\right)\ \dx x_d\\
&=\sqrt{\frac{n}{2\pi}}\exp(n\rho^2)\int_{\mathbb{R}}\exp\left(-\frac{n}{2}(x_d-2\rho)^2\right)\ \dx x_d\\
&=\exp(n\rho^2).
\end{align*}
Therefore inequality~\eqref{ChiLB} is satisfied (with equality) if the latter quantity is equal to $1+4(1-\eta)^2$, i.e. for
$$\rho= \sqrt{\frac{1}{n}\ln(1+4(1-\eta)^2)}.$$
This yields the claim.
\item Upper Bound.\\
Given $\alpha,\beta\in(0,\frac{1}{2})$, let $\delta=\min(\alpha,\beta)$ and $\tau_\delta=\sqrt{\frac{2}{n}v_\delta}$. Define the test
$$\varphi(X)=\mathds{1}_{\{X_d\geq \tau_\delta\}}.$$
Then for any $\mu\in\C$, we have
$$\PR_{\mu}(\varphi(X)=1)\leq \PR\left(\tfrac{1}{\sqrt{n}}\epsilon_d\geq \tau_\delta\right)\stackrel{(\ref{ConcIn}.\mathrm{I})}{\leq}\delta\leq\alpha.$$
On the other hand, let now $\rho=2\tau_\delta$. Then for any $\mu\in\mathcal{A}_\rho$, we have
$$\PR_{\mu}(\varphi(X)=0)\leq \PR\left(2\tau_\delta+\tfrac{1}{\sqrt{n}}\epsilon_d\leq \tau_\delta\right)\leq \PR\left(\tfrac{1}{\sqrt{n}}\epsilon_d\leq-\tau_\delta\right)\stackrel{(\ref{ConcIn}.\mathrm{I})}{\leq}\delta\leq\beta.$$
This concludes the proof since $\rho\eqsim\frac{1}{\sqrt{n}}$.
\end{enumerate}
\end{Proof}

\subsubsection{Proof of Theorem \ref{GenUpp}}\label{GenUppP}
\begin{Proof}
Given $\alpha,\beta\in(0,\frac{1}{2})$, let $\delta=\min(\alpha,\beta)$ and $\tau_\delta=\frac{d}{n}+\frac{2}{n}\sqrt{dv_\delta}+\frac{2}{n}v_\delta$. Define the test
$$\varphi(X)=\mathds{1}_{\{B(X,\sqrt{\tau_\delta})\cap\C=\varnothing\}}
=\mathds{1}_{\{\mathrm{dist}(X,\C)\geq \sqrt{\tau_\delta}\}}.$$
Then for any $\mu\in \C$, we have
$$\PR_{\mu}(\varphi(X)=1)\leq \PR \left( \Vert X - \mu \Vert \geq \sqrt{\tau_\delta} \right) \leq
\PR\left(\tfrac{1}{n}\Vert\epsilon\Vert^2\geq \tau_\delta\right)\stackrel{(\ref{ConcIn}.\mathrm{II})}{\leq}\delta\leq\alpha.$$
On the other hand, let now $\rho=2\sqrt{\tau_\delta}$ and $\mu\in\mathcal{A}_\rho$ arbitrary. Then analogously
$$\mathrm{dist}(X,\C)<\tau~\Rightarrow~\Vert X-\mu\Vert>\sqrt{\tau_{\delta}}$$ and hence
$$\PR_{\mu}(\varphi(X)=0)\leq \beta.$$
This concludes the proof since $\sqrt{\tau_\delta}\eqsim\sqrt{\frac{d}{n}}$.
\end{Proof}

\subsubsection{Proof of Theorem \ref{JNthm}}\label{JN}
\begin{Proof}
The arguments of this proof are related to the ones used in~\cite{Jud02} and~\cite{Cai11}. We decompose the proof into several steps.
\begin{enumerate}
\item Choice of priors.\\
We make use of the following lemma used and explained in \cite{Jud02} :\\
\begin{lemma0}
For any $M\in\mathbb{N}$ and $b>0$, there are distributions $\widetilde{\nu}_0$ and $\widetilde{\nu}_1$ with the following properties:
\begin{eqnarray}
&(\mathrm{I})&\supp(\widetilde{\nu}_0)\subseteq[-b,0],\ \ \supp(\widetilde{\nu}_1)\subseteq[-b,0]\cup\left\{\tfrac{b}{4M^2}\right\}\nonumber\\
&(\mathrm{II})&\widetilde{\nu}_1\left(\left\{\tfrac{b}{4M^2}\right\}\right)\geq\frac{1}{2} \label{JNPrior}\\
&(\mathrm{III})&\forall k\in\{0,1,\ldots,M\}:\ \int z^k\ \widetilde{\nu}_0(\dx z)=\int z^k\ \widetilde{\nu}_1(\dx z).\nonumber  
\end{eqnarray}
\end{lemma0}~\\
For now, let $\widetilde{\nu}_i$ be such distributions and $\nu_i=\widetilde{\nu}_i^{\otimes d}$, $i\in\{0,1\}$; $M,b$ and $\rho$ will be specified later. Furthermore, writing $\sigma^2=\frac{1}{n}$, let
$$\PR_i=\left(\widetilde{\nu}_i\ast\mathcal{N}(0,\sigma^2)\right)^{\otimes d},\ \ i\in\{0,1\},$$
where $\ast$ denotes convolution. Clearly, the corresponding density function can be written as
$$F_i(x)=\prod_{j=1}^d \left(\EW_{\mu_j\sim\widetilde{\nu}_i}[\phi(x_j;\mu_j,\sigma^2)]\right),\ \ i\in\{0,1\},$$
where $\phi(x;\mu,\sigma^2)$ is the density of $\mathcal{N}(\mu,\sigma^2)$. It will be convenient to examine the case $d=1$, denoted by $\widetilde{\PR}_i$.

\indent Note that $\PR_0$ is in accordance with $\PR_{\nu_0}$ from Section \ref{LBGen}, but the construction of $\nu_1$ does not warrant the notion of Euclidean distance we are interested in - $\nu_1$ has support inside $\C_{\mathrm{O}}$ - hence the slight difference in notation. This technical obstacle is necessary for the property (\ref{JNPrior}.III), but it can be resolved for a small price, which we explain in the last step of this proof.
\item Controlling the total variation distance.\\
Based on our construction, we have for $i\in\{0,1\}$ and fixed $x\in\mathbb{R}$
\begin{align}\EW_{\mu\sim\widetilde{\nu}_i}[\phi(x;\mu,\sigma^2)]&=\phi(x;0,\sigma^2)\int \exp\left(\frac{2x\mu-\mu^2}{2\sigma^2}\right)\widetilde{\nu}_i(\dx \mu) \nonumber\\ 
&=\phi(x;0,\sigma^2)\int \sum_{k=0}^\infty \frac{1}{k!}\left(\frac{2x\mu-\mu^2}{2\sigma^2}\right)^k\ \widetilde{\nu}_i(\dx \mu)\nonumber\\ 
&=\phi(x;0,\sigma^2)\int \sum_{k=0}^\infty \frac{1}{k!(2\sigma^2)^k}(2x\mu-\mu^2)^{k}\ \widetilde{\nu}_i(\dx \mu).\label{DensBound}
\end{align}
Let now
$$D_k(x):=\int (2x\mu-\mu^2)^{k}\ \widetilde{\nu}_1(\dx \mu)-\int(2x\mu-\mu^2)^{k}\ \widetilde{\nu}_0(\dx \mu).$$
Then \eqref{DensBound} in conjunction with (\ref{JNPrior}.III) tells us that
$$\frac{\EW_{\mu\sim\widetilde{\nu}_1}[\phi(x;\mu,\sigma^2)]-\EW_{\mu\sim\widetilde{\nu}_0}[\phi(x;\mu,\sigma^2)]}{\phi(x;0,\sigma^2)}=\sum_{k=\lfloor M/2\rfloor+1}^\infty \frac{1}{k!(2\sigma^2)^k}D_k(x).$$
and thus
\begin{align}
	\Vert \widetilde{\PR}_1-\widetilde{\PR}_0\Vert_{\mathrm{TV}}&=\int\left|\EW_{\mu\sim\widetilde{\nu}_1}[\phi(x;\mu,\sigma^2)]-\EW_{\mu\sim\widetilde{\nu}_0}[\phi(x;\mu,\sigma^2)]\right|\ \dx x \nonumber\\
	&\leq \sum_{k=\lfloor M/2\rfloor+1}^\infty \frac{1}{k!(2\sigma^2)^k}\left|\int D_k(x)\phi(x;0,\sigma^2)~\dx x\right|\label{bound1}
	\end{align}
We take a moment to upper bound the individual summands: Since
$$(2x\mu-\mu^2)^{k}\leq 4^{k}|x|^kb^k+2^k b^{2k}$$
and, with \cite{Win12}, 
$$\int |x|^k \phi(x,0,\sigma^2)\ \dx x= \frac{\sigma^k \sqrt{2}^k}{\sqrt{\pi}}\Gamma((k+1)/2)\leq \frac{\sigma^k \sqrt{2}^k}{\sqrt{\pi}}\left\lceil\frac{k}{2}\right\rceil !,$$
we have
\begin{align*}\left|\int D_k(x)\phi(x;0,\sigma^2)~\dx x\right|&\leq 2\left(4^kb^k\int|x|^k\phi(x;0,\sigma^2)~\dx x+2^kb^{2k}\right)\\
&\leq 2\left(\frac{1}{\sqrt{\pi}}\left(4\sqrt{2}b\sigma\right)^k\left\lceil\frac{k}{2}\right\rceil !+2^kb^{2k}\right).\end{align*}
Now through Stirling's approximation and elementary manipulation, with $M\geq 32$ we obtain
\begin{align*}
\displaystyle\dfrac{\left\lceil\frac{k}{2}\right\rceil !}{k!}&\leq \dfrac{\mathrm{e}}{\sqrt{2\pi}}\dfrac{\left\lceil\frac{k}{2}\right\rceil^{\left\lceil\frac{k}{2}\right\rceil}}{k^{k+\frac{1}{2}}}\underbrace{\frac{\left\lceil\frac{k}{2}\right\rceil^{1/2}}{\mathrm{e}^{\left\lceil\frac{k}{2}\right\rceil}}}_{\leq \sqrt{3}/k}\mathrm{e}^k\\
&\leq \dfrac{\mathrm{e}\sqrt{3}}{\sqrt{2\pi}}\left(\sqrt{\frac{k+1}{2k^2}}\right)^{k+1}\mathrm{e}^k\\
&\leq \dfrac{\mathrm{e}\sqrt{3}}{5\sqrt{2\pi}}\left(\sqrt{\frac{17}{32}}\frac{1}{\sqrt{k}}\right)^k\\
&\leq \frac{1}{2}\left(\sqrt{\frac{17}{32}}\frac{1}{\sqrt{k}}\right)^k
\end{align*}
and
$$k!\geq \sqrt{2\pi}k^k\sqrt{k}\mathrm{e}^{-k}\geq 4\sqrt{\pi}\left(\frac{k}{\mathrm{e}}\right)^k.$$
That yields
\begin{align*}\frac{1}{k!(2\sigma^2)^k}\left|\int D_k(x)\phi(x;0,\sigma^2)~\dx x\right|&\leq 2\left(\frac{1}{k!\sqrt{\pi}}\left(\frac{4}{\sqrt{2}}\frac{b}{\sigma}\right)^k\left\lceil\frac{k}{2}\right\rceil !+\frac{1}{k!}\left(\frac{b}{\sigma}\right)^{2k}\right)\\
&\leq 2\left(\frac{1}{2}\left(\frac{\sqrt{17}}{2}\frac{b}{\sigma\sqrt{k}}\right)^k+\frac{1}{4\sqrt{\pi}}\left(\mathrm{e}\frac{b^2}{\sigma^2k}\right)^k\right)\end{align*}
At this point, we introduce a more explicit choice of $b$, namely $b=c\sqrt{M}\sigma$ with $c=\frac{2\sqrt{2}}{\sqrt{17}\mathrm{e}}\geq\frac{1}{4}$. This choice guarantees
\begin{align*}\frac{1}{k!(2\sigma^2)^k}\left|\int D_k(x)\phi(x;0,\sigma^2)~\dx x\right|
&\leq \left(1+\frac{1}{2\sqrt{\pi}}\right)\left(\frac{\sqrt{17}}{2}\frac{b}{\sigma\sqrt{k}}\right)^k\end{align*}
and moreover, continuing \eqref{bound1}, 
\begin{align*}
	\Vert \widetilde{\PR}_1-\widetilde{\PR}_0\Vert_{\mathrm{TV}}&\leq  \left(1+\frac{1}{2\sqrt{\pi}}\right)\sum_{k=\lfloor M/2\rfloor+1}^\infty\left(\frac{\sqrt{17}}{2}\frac{b}{\sigma\sqrt{k}}\right)^k\\
	&\leq \left(1+\frac{1}{2\sqrt{\pi}}\right)\frac{2}{\mathrm{e}-2}\left(\frac{2}{\mathrm{e}}\right)^{\lfloor M/2\rfloor}
\end{align*}
and hence finally
$$\Vert \PR_1-\PR_0\Vert_{\mathrm{TV}}\leq d\left(1+\frac{1}{2\sqrt{\pi}}\right)\frac{2}{\mathrm{e}-2}\left(\frac{2}{\mathrm{e}}\right)^{\lfloor M/2\rfloor}.$$
By direct computation, we now see that for any $\eta'\in(0,1)$
$$\frac{1}{2}\Vert \PR_1-\PR_0\Vert_{\mathrm{TV}}\leq 1-\eta'$$
is fulfilled if
$$M\geq \frac{2}{1-\ln(2)}\ln(d)+\underbrace{1+\frac{2}{1-\ln(2)}\ln\left(\frac{1.8}{1-\eta'}\right)}_{=:c_{\eta'}},$$
so we choose
$$M:=\max\left(32,\left\lceil \frac{2}{1-\ln(2)}\ln(d)+c_{\eta'}\right\rceil\right).$$
\item Application.\\
Note that this upper bound $\frac{1}{2}\Vert \PR_1-\PR_0\Vert_{\mathrm{TV}}\leq1-\eta'$ does not formally allow for determining a lower bound on $\rho^\ast(\C)$ yet since $H_0$ and $H_1$ are not separated in a Euclidean sense. In a final step, we will resolve this by a suitable restriction of $H_1$.\\

Let $Y=\sum_{i=1}^d \mathds{1}_{\{\mu_i=u\}}$, i.e. the number of coordinates of $\mu$ taking the value $u=\frac{b}{4M^2}$. Obviously, if $\mu\sim\nu_1$, we have $Y\sim\mathrm{Bin}(d,\widetilde{\nu}_1(\{u\}))$. By property (\ref{JNPrior}.II) and Hoeffding's inequality, this yields that if $d\geq 42$, 
$$\PR_{\mu\sim \nu_1}\left(Y\geq\tfrac{d}{3}\right)\geq \frac{9}{10}.$$ Now, let $\xi=\{Y\geq\frac{d}{3}\}$ and 
$$H_1':\ \mu\sim \nu_1|\xi.$$
Assuming that for some test $\varphi$ the relation
$$\PR_{\mu\sim\nu_0}(\varphi=1)+\PR_{\mu\sim\nu_1}(\varphi=0)\geq \eta'$$
holds, we can conclude
\begin{eqnarray*}
\PR_{\nu_0}(\varphi=1)+\PR_{\nu_1|\xi}(\varphi=0)&=&\PR_{\nu_0}(\varphi=1)+1-\frac{\PR_{\nu_1}(\{\varphi=1\}\ \cap\ \xi)}{\PR_{\nu_1}(\xi)}\\
&\geq &\PR_{\nu_0}(\varphi=1)+1-\frac{10}{9}\PR_{\nu_1}(\varphi=1)\\
&\geq &\PR_{\nu_0}(\varphi=1)+\frac{10}{9}\PR_{\nu_1}(\varphi=0)-\frac{1}{9}\\
&\geq&\PR_{\nu_0}(\varphi=1)+\PR_{\nu_1}(\varphi=0)-\frac{1}{9}\\
&\geq&\eta'-\frac{1}{9}.
\end{eqnarray*}
Hence, inference from the testing problem discussed in Steps 1 and 2 to the problem $$H_0:\ \mu\sim\nu_0\ \ \mathrm{vs.}\ \ H_1':\ \mu\sim\nu_1|\xi$$
is valid as long as $\eta\in(0,\frac{8}{9})$ ($\eta$ corresponds to $\eta'-\frac{1}{9}$ above).\\

The following observation concludes the proof: Clearly, $H_1'$ agrees with $\mathcal{A}_\rho$ for 
$$\rho=\frac{\sqrt{d}}{\sqrt{3}}\frac{b}{4M^2}=\frac{1}{2\mathrm{e}\sqrt{25.5}}\frac{1}{M^{3/2}}\frac{\sqrt{d}}{\sqrt{n}}´\geq\frac{1}{28}\frac{1}{M^{3/2}}\frac{\sqrt{d}}{\sqrt{n}}.$$
\end{enumerate}
\end{Proof}

\subsection{Proofs for Section \ref{Smooth}}

\subsubsection{Proof of Lemma \ref{Tay}}
\begin{Proof}
We need to ensure that on $B$, the graph of $f$ remains below $\widetilde{B}=B_d([\mathds{O}_{d-1},R],R)$ since that corresponds to the fact that $\widetilde{B}$ is locally contained in $\C$, as required in Definition \ref{Rrounded}. This is equivalent to
$$\forall x\in B:~0\leq f(x)\leq R-\sqrt{R^2-\Vert x\Vert^2}.$$
Applying Taylor's theorem with Lagrange's remainder yields
$$\exists s\in(0,1):\ f(x)=\frac{1}{2} x^T Hf(sx) x,$$
since by construction $f(\mathds{O}_{d-1})=0$ and $\nabla f(\mathds{O}_{d-1})=\mathds{O}_{d-1}$. Clearly, in order that $f\geq 0$ on $B$, it is sufficient to require $\lambda_{\mathrm{min}}(Hf(y))\geq 0$ for $y\in B\backslash\{\mathds{O}_{d-1}\}$. On the other hand, we can use a classical eigenvalue representation to obtain the desired upper bound: For some $s\in(0,1)$,
\begin{eqnarray*}f(x)&=& \frac{1}{2}\Vert x\Vert^2 \left(\frac{x}{\Vert x\Vert}\right)^T Hf(sx)\left(\frac{x}{\Vert x\Vert}\right)\\
&\leq&\frac{1}{2}\Vert x\Vert^2\max_{\Vert y\Vert=1} y^THf(sx)y\\
&=&\frac{1}{2}\Vert x\Vert^2 \lambda_{\mbox{max}}(Hf(sx))\\
&\leq&\frac{1}{2R}\Vert x\Vert^2\\
&\leq&R-\sqrt{R^2-\Vert x\Vert^2}\end{eqnarray*}
by assumption and \eqref{SqrtIn2}.
\end{Proof}
\subsubsection{Proof of Theorem \ref{GenRSm}}
\begin{Proof}
We define the test statistic
$$T(X):=\dist(X,\C),$$
and a corresponding test of the form $\varphi(X)=\mathds{1}_{\{T(X)\geq\tau\}}$.\\
Let $\mu\in\C$. W.l.o.g. assume that $\mu':=\mathds{O}_d\in\C\subseteq\mathbb{R}^{d-1}\times[0,\infty)$ and $\mu'$ minimises the distance between $\mu$ and $\partial\C$. Now let $z=[\mathds{O}_{d-1},R]$ so that by construction $\mu'\in \overline{B_d(z,R)}\subseteq\C$. 

For $\tau>0$, we have
\begin{eqnarray*}\mbox{dist}(X,\C)\geq \tau &\ \Longrightarrow\ &\mbox{dist}(X,B_d(z,R))\geq \tau\\
&\ \Longrightarrow\ &\left\Vert \tfrac{1}{\sqrt{n}}\epsilon-z\right\Vert-R\geq \tau.\\
\end{eqnarray*}
Now, writing $\epsilon_{1:(d-1)}:=[\epsilon_1,\epsilon_2,\ldots,\epsilon_{d-1}]$ and  using (\ref{SqrtIn1}), we obtain
\begin{eqnarray*}
\left\Vert \tfrac{1}{\sqrt{n}}\epsilon-z\right\Vert-R&=&\sqrt{\left\Vert \tfrac{1}{\sqrt{n}}\epsilon-z\right\Vert^2}-R\\
&=&\sqrt{\tfrac{1}{n}\Vert \epsilon_{1:(d-1)}\Vert^2+\left(\tfrac{1}{\sqrt{n}}\epsilon_{d}-R\right)^2}-R\\
&\leq&\sqrt{\tfrac{1}{n}\Vert \epsilon_{1:(d-1)}\Vert^2+\left(R+\tfrac{1}{\sqrt{n}}|\epsilon_{d}|\right)^2}-R\\
&\leq&R+\tfrac{1}{\sqrt{n}}|\epsilon_{d}|+\tfrac{\Vert \epsilon_{1:(d-1)}\Vert^2}{2n(R+\tfrac{1}{\sqrt{n}}|\epsilon_{d}|)}-R\\
&\leq&\tfrac{1}{\sqrt{n}}|\epsilon_{d}|+\tfrac{\Vert \epsilon_{1:(d-1)}\Vert^2}{2nR},\\
\end{eqnarray*}
which tells us
$$\PR_\mu(T(X)\geq\tau)\leq \PR\left(\tfrac{1}{\sqrt{n}}|\epsilon_{d}|+\tfrac{\Vert \epsilon_{1:(d-1)}\Vert^2}{2nR}\geq\tau\right).$$
Clearly, this bound holds generally in the sense
$$\sup\limits_{\mu\in\C}\PR_\mu(\varphi(X)=1)\leq \PR\left(\tfrac{1}{\sqrt{n}}|\epsilon_{d}|+\tfrac{\Vert \epsilon_{1:(d-1)}\Vert^2}{2nR}\geq\tau\right).$$\\
Based on the general property 
$$\PR(A\geq \tau_1)\leq\tfrac{\alpha}{2}\ \wedge\ \PR(B\geq\tau_2)\leq\tfrac{\alpha}{2}\  \ \Longrightarrow\ \ \PR(A+B\geq\tau_1+\tau_2)\leq\alpha$$
for random variables $A$ and $B$ and by using (\ref{ConcIn}.I) and (\ref{ConcIn}.II), we finally obtain the rejection threshold
$$\tau:=\sqrt{\tfrac{2}{n}v_{\alpha/4}}+\tfrac{d}{2nR}+\tfrac{2}{nR}\sqrt{dv_{\alpha/2}}+\tfrac{1}{nR}v_{\alpha/2}\eqsim\max\left(\frac{1}{\sqrt{n}},\frac{d}{nR}\right)$$
for a fixed level $\alpha\in(0,\frac{1}{2})$.

On the other hand, w.l.o.g., choose $\mu=[\mathds{O}_{d-1},-\rho]$. This is valid since by construction $\mu$ minimises the distance between $\C$ and $\mathcal{A}_\rho$ and $\mathds{O}_d$ represents an arbitrary element of $\partial\C$. We have
$$\dist(X,\C)\leq\tau\ \Longrightarrow\ X_d\geq -\tau\ \Longleftrightarrow\ \epsilon_d\geq \sqrt{n}(\rho-\tau),$$
so that it is sufficient to ensure
$$\sup_{\mu\in\mathcal{A}_\rho}\PR_\mu(\varphi(X)=0)\ \leq\ \PR(\epsilon_d\geq\sqrt{n}(\rho-\tau))\ \leq\ \beta\in(0,\tfrac{1}{2}),$$
which leads to the condition
$$\rho\geq \tau+\sqrt{\tfrac{2}{n}v_{\beta}}\sim\tau.$$
This concludes the proof.
\end{Proof}

\subsubsection{Proof of Theorem \ref{IO}}
\begin{Proof}
This is a variation on the proof of Theorem~\ref{JNthm}. Using the same construction
and notation as previously, and taking $d\geq 3$, let now for $i\in\{0,1\}$
$$\nu_i=\widetilde{\nu}_i^{\otimes d-1}\otimes \delta_R.$$
Since the mutual deterministic coordinate $\mu_d=R$ is irrelevant for the total variation distance between the resulting distributions $\PR_0$ and $\PR_1$, the bounds in Step 2 of the proof of Theorem~\ref{JNthm} also hold here with $d-1$ instead of $d$.\\
The most important modification arises when calculating $\rho$: Now, if at least $\frac{d-1}{3}$ of the coordinates take the value $u=\frac{b}{4M^2}$, computing the Euclidean distance of $\mu$ from $\mathcal{C}$ and using (\ref{SqrtIn1}) leads to
\begin{eqnarray*}\rho^\ast(\C)\geq\sqrt{R^2+\frac{d-1}{3}u^2}-R&\geq& \frac{(d-1)u^2}{6\sqrt{R^2+\frac{d-1}{3}u^2}}\\
&\geq& \frac{(d-1)u^2}{6R+\frac{6}{\sqrt{3}}\sqrt{d-1}u}\\
&\geq&\frac{1}{12}\min\left(\frac{(d-1)u^2}{R},\sqrt{3}\sqrt{d-1}u\right)\\
&\sim&\min\left(\frac{1}{\ln(d)^3}\cdot \frac{d}{nR},\frac{1}{\ln(d)^{3/2}}\cdot \frac{\sqrt{d}}{\sqrt{n}}\right),\end{eqnarray*}
if $d$ is large enough in the sense that $M_{\eta}\leq C\ln(d-1)$ for some $C>0$, where $M_{\eta}$ is given in the statement of the theorem.
This concludes the proof.
\end{Proof}

\subsection{Proofs for Section \ref{Dis}}
\subsubsection{Proof of Theorem \ref{Ball}}
\begin{Proof}
W.l.o.g., let $z=\mathds{O}_d$. We prove independently that $\rho^\ast(\C)$ is lower and upper bounded by the right hand side of (\ref{BallRate}).

\paragraph{Lower Bound.} Let $\nu_0=\delta_{Re_d}$, giving rise to the density function
$$F_{\nu_0}(x):=\left(\frac{n}{2\pi}\right)^{\frac{d}{2}}\exp\left(-\frac{n}{2}(x_d-R)^2\right)\prod_{i=1}^{d-1} \exp\left(-\frac{n}{2}x_i^2\right).$$

On the other hand, for a suitable $h>0$ specified in a moment, let $\nu_\rho$ be the uniform distribution on
$$P_h:=\{[h\cdot v,R]\ |\ v\in\{-1,1\}^{d-1}\}$$
Since each element of $P_h$ has Euclidean distance $\sqrt{R^2+(d-1)h^2}-R$ from $\C$, which should correspond to $\rho$, we set $h^2=\frac{(R+\rho)^2-R^2}{d-1}$. This gives rise to the following density function:
\begin{eqnarray*}
F_{\nu_\rho}(x)&:=&\left(\frac{n}{2\pi}\right)^{\frac{d}{2}}\exp\left(-\frac{n}{2}(x_d-R)^2\right)\frac{1}{2^{d-1}}\sum_{v_1,\ldots,v_{d-1}\in\{-1,1\}}\prod_{i=1}^{d-1} \exp\left(-\frac{n}{2}(x_i-h\cdot v_i)^2\right)\\
&=&\left(\frac{n}{2\pi}\right)^{\frac{d}{2}}\exp\left(-\frac{n}{2}(x_d-R)^2\right)\frac{1}{2^{d-1}}\prod_{i=1}^{d-1}\exp\left(-\frac{n}{2}x_i^2-\frac{n}{2}h^2\right)2\cosh\left(nhx_i\right)\\
&=&\left(\frac{n}{2\pi}\right)^{\frac{d}{2}}\exp\left(-\frac{n}{2}(x_d-R)^2-(d-1)\frac{n}{2}h^2\right)\prod_{i=1}^{d-1}\exp\left(-\frac{n}{2}x_i^2\right)\cosh\left(nhx_i\right),
\end{eqnarray*}
so that
$$F_{\nu_\rho}^2(x):=\left(\frac{n}{2\pi}\right)^{d}\exp\left(-n(x_d-R)^2-(d-1)nh^2\right)\prod_{i=1}^{d-1}\exp\left(-nx_i^2\right)\cosh^2\left(nhx_i\right).$$
Now, using the fact that $\EW[\cosh^2(aY)]=\exp(a^2\sigma^2)\cosh(a^2\sigma^2)$ for $Y\sim\mathcal{N}(0,\sigma^2)$, we have
\begin{eqnarray*}
\int_{\mathbb{R}^d}\tfrac{F_{\nu_\rho}^2(x)}{F_{\nu_0}(x)}\ \dx x&=&\int_{\mathbb{R}^d}\left(\tfrac{n}{2\pi}\right)^{\frac{d}{2}}\exp\left(-\tfrac{n}{2}(x_d-R)^2-(d-1)nh^2\right)\prod_{i=1}^{d-1}\cosh^2\left(nhx_i\right)\exp\left(-\tfrac{n}{2}x_i^2\right)\ \dx x\\
&=&\left(\tfrac{n}{2\pi}\right)^{\frac{d}{2}}\exp\left(-(d-1)nh^2\right)\int_{\mathbb{R}}\exp\left(-\tfrac{n}{2}(x_d-R)^2\right)\dx x_1\\
&&\ \ \ \ \ \ \ \ \ \ \ \ \ \ \ \ \ \ \ \ \ \ \ \ \ \ \ \ \ \ \ \ \ \ \ \ \ \ \ \cdot\left[ \int_{\mathbb{R}}\cosh^2\left(nhx\right)\exp\left(-\tfrac{n}{2}x^2\right)\dx x\right]^{d-1}\\
&=&\exp\left(-(d-1)nh^2\right)\cdot\left[\exp(nh^2)\cosh(nh^2)\right]^{d-1}\\
&=&\cosh(nh^2)^{d-1}.
\end{eqnarray*}
Now, by Taylor expansion we obtain the bound
\begin{equation}nh^2\leq 1\ \ \Longrightarrow\ \ \cosh(nh^2)\leq 1+\frac{e}{2}n^2h^4,\label{Tay1}\end{equation}
so that
\begin{equation}
\label{sidebound}
\ln\left(\cosh(nh^2)^{d-1}\right)\leq (d-1)\frac{e}{2}n^2h^4=\frac{n^2}{(d-1)}\frac{e}{2}((R+\rho)^2-R^2)^2,
\end{equation}
whenever
\begin{equation}
\label{sidecond}
nh^2\leq 1\ \ \mathrm{i.e.}\ \ \rho\leq\sqrt{\frac{d-1}{n}+R^2}-R.
\end{equation}
The upper bound \eqref{sidebound} leads to the following condition: 
\begin{equation}\rho\leq \sqrt{\frac{\sqrt{d-1}}{n}s+R^2}-R,
\text{ where }  s:= \sqrt{\frac{2}{e}\ln(1+4(1-\eta)^2)}, \label{Tay2}\end{equation}
which is sufficient for ensuring (\ref{ChiLB}), provided that \eqref{sidecond} holds. It is straightforward to see that (\ref{Tay2}) implies \eqref{sidecond} as long as 
$$d\geq 1+\frac{2}{e}\ln(5)\ \ \mathrm{i.e.}\ \ d\geq 3.$$
It remains to investigate $(\ref{Tay2})$ a little closer. Application of (\ref{SqrtIn1}) now yields the following, defining $t>0$ via the relation $R^2 = t^2 s\frac{\sqrt{d-1}}{n}$:
\begin{align*}
\sqrt{\frac{\sqrt{d-1}}{n}s+R^2}-R\geq \frac{\sqrt{s}}{2\sqrt{1+t^2}}\frac{(d-1)^{\frac{1}{4}}}{n^{\frac{1}{2}}}&\eqsim \frac{(d-1)^{\frac{1}{4}}}{n^{\frac{1}{2}}} \min(t^{-1},1) \\
&\eqsim \min\left(\frac{(d-1)^{\frac{1}{4}}}{n^{\frac{1}{2}}},\frac{(d-1)^{\frac{1}{2}}}{n\cdot R}\right).
\end{align*}
If on the one hand $R\lesssim \frac{(d-1)^{\frac{1}{4}}}{n^{\frac{1}{2}}}\sqrt{s}$, that is  $R \leq  t\sqrt{s}\frac{(d-1)^{\frac{1}{4}}}{n^{\frac{1}{2}}}$for $t>0$, we have
$$\sqrt{\frac{\sqrt{d-1}}{n}s+R^2}-R\geq \frac{\sqrt{s}}{2\sqrt{1+t^2}}\frac{(d-1)^{\frac{1}{4}}}{n^{\frac{1}{2}}}\eqsim \min\left(\frac{(d-1)^{\frac{1}{4}}}{n^{\frac{1}{2}}},\frac{(d-1)^{\frac{1}{2}}}{n\cdot R}\right).$$
Analogously, the case $R\gtrsim \frac{(d-1)^{\frac{1}{4}}}{n^{\frac{1}{2}}}\sqrt{s}$ also yields
$$\sqrt{\frac{\sqrt{d-1}}{n}s+R^2}-R\gtrsim \min\left(\frac{(d-1)^{\frac{1}{4}}}{n^{\frac{1}{2}}},\frac{(d-1)^{\frac{1}{2}}}{n\cdot R}\right).$$

\paragraph{Upper Bound.} We define the test statistic
$$T(X):=\Vert X\Vert^2-R^2$$
and a corresponding test of the form $\varphi(X)=\mathds{1}_{\{T(X)\geq\tau\}}$.
On the one hand, in order to control the type-I-error probability, take any $\mu\in C$, so that $\Vert \mu\Vert\in[0,R]$. Clearly, $n\Vert X\Vert^2\sim\chi^2_{n\Vert\mu\Vert^2}(d)$. Therefore, for $\tau'>0$ and with the notation $Z_{\lambda}\sim \chi^2_{\lambda}(d)$, we can guarantee
$$\PR_{\mu}(\varphi(X)=1)=\PR(Z_{n\Vert\mu\Vert^2}\geq n\Vert\mu\Vert^2+n\tau')\leq\alpha$$
by setting
$$\tau'=\frac{d}{n}+2\sqrt{\left(\frac{d}{n^2}+\frac{2}{n}\Vert\mu\Vert^2\right)v_{\alpha}}+\frac{2}{n}v_{\alpha},$$
where we use (\ref{ConcIn}.II). Since $\Vert\mu\Vert\leq R$, this yields that
$$\sup_{\mu\in\C}\PR_{\mu}(\varphi(X)=1)\leq \sup_{\mu\in\C}\PR(Z_{n\Vert\mu\Vert^2}\geq nR^2+n\tau)\leq\alpha$$
for
$$\tau=\frac{d}{n}+2\sqrt{\left(\frac{d}{n^2}+\frac{2}{n}R^2\right)v_{\alpha}}+\frac{2}{n}v_{\alpha}.$$
On the other hand, in order to satisfy a prescribed level $\beta\in(0,\frac{1}{2})$ for the Type-II-error, take any $\mu\in\mathcal{A}_\rho$ with $\Vert\mu\Vert\geq R+\rho$. Then again, $n\Vert X\Vert^2\sim\chi^2_{n\Vert\mu\Vert^2}(d)$, so that we need to ensure
\begin{equation}\PR_\mu(\varphi(X)=0)=\PR(Z'\leq nR^2+n\tau)\leq\beta,\ \ \mathrm{where\ }Z'\sim\chi^2_{n\Vert\mu\Vert^2}(d).\label{BallTII}\end{equation}
In this case, (\ref{ConcIn}.III) yields the sufficient condition
$$d+nR^2+2\sqrt{(d+2nR^2)v_{\alpha}}+2v_\alpha \leq d+n\Vert\mu\Vert^2-2\sqrt{(d+2n\Vert\mu\Vert^2)v_{\beta}}.$$
The right hand side is increasing in $\Vert\mu\Vert$ if $d\geq\ln(1/\beta)$, so that, similar as for the type-I-error, \eqref{BallTII} holds uniformly over $\mathcal{A}_{\rho}$ if
$$d+nR^2+2\sqrt{(d+2nR^2)v_{\alpha}}+2v_\alpha \leq d+n(R+\rho)^2-2\sqrt{(d+2n(R+\rho)^2)v_{\beta}}.$$
Using $\sqrt{a+b}\leq \sqrt{a}+\sqrt{b}$ (for $a,b>0$) and (\ref{SqrtIn1}) respectively, we obtain two different sufficient bounds for $\rho$:
\begin{eqnarray*}
\rho&\geq&\sqrt{2}\frac{d^{\frac{1}{4}}}{\sqrt{n}}(\sqrt{v_\alpha}+\sqrt{v_\beta})+\sqrt{\frac{2}{n}}(\sqrt{v_\alpha}+2\sqrt{v_\beta})\,;\\
\rho&\geq&\frac{\sqrt{d}}{nR+2\sqrt{nv_\alpha}}(\sqrt{v_\alpha}+\sqrt{v_\beta})+\sqrt{\frac{2}{n}}(\sqrt{v_\alpha}+2\sqrt{v_\beta}).
\end{eqnarray*}
Therefore, as claimed, the upper bound
$$\rho^\ast(\C)\lesssim\max\left(\frac{1}{\sqrt{n}},\min\left(\frac{d^{\frac{1}{4}}}{\sqrt{n}},\frac{\sqrt{d}}{nR}\right)\right)$$
holds. This concludes the proof.
\end{Proof}

\bibliographystyle{acm}
\bibliography{Biblio}
%

\end{document}